\documentclass[12pt,a4paper,twoside]{amsart}
\usepackage{amsfonts, amsthm, amsmath, amssymb}
\usepackage{mathrsfs,amsmath}
\usepackage{hyperref}
\usepackage{bbm}
\hypersetup{colorlinks=false}

\usepackage[margin=2.2cm]{geometry}

\usepackage{helvet}
\usepackage[thinc]{esdiff}

\usepackage{amsthm}
\usepackage{lineno}
\newtheorem{theorem}{Theorem}
\newtheorem{lemma}{Lemma}[subsection]
\newtheorem{remark}{Remark}

\newtheorem{proposition}{Proposition}[subsection]
\newtheorem{corollary}{Corollary}

\usepackage{accents}
\newlength{\dhatheight}
\newcommand{\doublehat}[1]{%
    \settoheight{\dhatheight}{\ensuremath{\hat{#1}}}%
    \addtolength{\dhatheight}{-0.35ex}%
    \hat{\vphantom{\rule{1pt}{\dhatheight}}%
    \smash{\hat{#1}}}}

\begin{document}

\author{ Aritra Ghosh }
\title{Weight aspect asymptotic and simultaneous non-vanishing for Rankin-Selberg $L$-functions }
\address{Aritra Ghosh \newline  TAN group, École Polytechnique Fédérale de Lausanne (EPFL), MA C3 634, Bâtiment MA, Station 8, 1015 Lausanne, Switzerland; email: aritrajp30@gmail.com}

\begin{abstract}
In this article we show simultaneous non-vanishing of two Rankin-Selberg $L$-functions by proving an asymptotic result in weight aspect. The main input of this paper is to remove the $t$-integral dependence from the result of Blomer-Harcos (see \cite{BH2}) and getting a square root exponent for the error term.
\end{abstract}

\maketitle
\tableofcontents

\section{Introduction}
 An interesting and central problem of number theory is to prove vanishing or non-vanishing properties of an automorphic $L$-function which also encode a deep arithmetic property, seen for many cases. For example, one can see that the Birch and Swinnerton-Dyer conjecture relates the order of vanishing of the Hasse-Weil $L$-function at the central point to the rank of an elliptic curve.

 If $\pi$ is taken as a unitary cuspidal automorphic representation of $\mathrm{GL(n)}$ with $s\in \mathbb{C}$  then the non-vanishing of $L(s,\pi)$ twists by a Dirichlet character $\chi$ had been studied many mathematicians (see \cite{BR}, \cite{JN}, \cite{Luo}, \cite{LRS},  \cite{Ro}). Specially the value of $L$-function at the central point $s= \frac{1}{2}$ are the subject of intensive studies in various aspects. Also, positivity of the Dirichlet $L$-functions with real characters at $s = \frac{1}{2}$ would give remarkable lower bounds for the class number of imaginary quadratic fields. Also one can see that, a good positive lower bound for the central values of Hecke $L$-functions would deny the existence of Landau-Siegel zero. For other non-vanishing results one can see \cite{BHM}, \cite{DK}, \cite{DFI2}, \cite{AGL}, \cite{HM}, \cite{HLX}, \cite{RK}, \cite{KMV}, \cite{XL}, \cite{SCL}, \cite{GK}, \cite{PGM}, \cite{MV}, \cite{MS}, \cite{RR}, \cite{ZX}. In the current development of the generalized Ramanujan conjecture (see \cite{LRS}), the non-vanishing of certain Rankin-Selberg $L$-functions plays a crucial role. 

 Let $H_{k}$ denote the set of holomorphic Hecke cusp forms $g$ of weight $k$ for $\mathrm{SL}(2, \mathbb{Z})$, where $g(z)$ has the normalized (i.e., $\lambda_{g}(1)=1$) Fourier expansion given below:

 $$g(z)=\sum_{n=1}^{\infty} \lambda_{g}(n)n^{\frac{k-1}{2}}e^{2\pi i n z}.$$

\noindent
For a cusp form $g$, we take 
$$\omega_{g}:=\frac{\Gamma(k-1)}{(4\pi)^{k-1} \|g\|^2 },$$

\noindent
where
$$\|g\|^2 := \int_{\mathrm{SL}(2,\mathbb{Z})\backslash \mathbb{H}}\lvert g(z)\rvert^2 y^{k-2} \; \mathrm{d}x \; \mathrm{d}y,$$

\noindent
is the usual Petersson norm. For the Hecke eigenform $g$, it is known by the work of Hoffstein–Lockhart and Iwaniec \cite{HL, Iws} that

$$\omega_{g}\asymp k^{-1+\epsilon}.$$

Let $f,h$ be two distinct fixed automorphic forms for $\mathrm{SL}(2,\mathbb{Z})$ and $p$ be an odd prime number. Then the $L$-function associated with $f \otimes g $ is given by
$$L(s, f \otimes g ) =  \zeta_{p}(2s)\sum_{n=1}^{\infty}\frac{\lambda_{f}(n)\lambda_{g}(n)}{n^s}, \, \, \, \, \zeta_{p}(s)= \prod_{q\neq p }\left( 1-\frac{1}{q^s}\right)^{-1},$$
for $\mathrm{Re}(s) > 1$, which can be analytically extended to the whole complex plane $\mathbb{C}$ and satisfies a functional equation relating $s$ and $1-s$ and for the last product in the previous line, $q$ runs over all primes except $p$. Thousgh for us every form has full level.

\noindent
Here in this article we consider the sum
$$\mathcal{I}(K):=\sum_{k \equiv 0 \bmod 2}u\left(\frac{k-1}{K} \right)\sum_{g \in H_{k}}\omega_{g}^{-1}L(\frac{1}{2},f \otimes g)L(\frac{1}{2},g \otimes h)  .$$

\noindent
Our main theorem is the following asymptotic result:
\begin{theorem}\label{main}
    Let $f,h$ be two distinct fixed automorphic forms for $\mathrm{SL}(2,\mathbb{Z})$. Let $u \in \mathcal{C}_{c}^{\infty}(0,\infty)$. Then for $K$ large, we have
    \begin{equation*}
        \begin{split}
            & \mathcal{I}(K):=\sum_{k \equiv 0 \bmod 2}u\left(\frac{k-1}{K} \right)\sum_{g \in H_{k}}\omega_{g}^{-1}L(\frac{1}{2},f \otimes g)L(\frac{1}{2},g \otimes h) \\
             & = KP(\log K )  + \mathrm{O}(K^{\frac{1}{2}+\epsilon}).
             \end{split}
             \end{equation*}
             \noindent
    where $P(\log x)$ is a polynomial in $\log x$, depending on $f,h$ only. $P(\log x)$ is a degree $2$ polynomial if $h\neq \overline{f}$, atleast one of $f,h$ is a cusp form; $P(\log x)$ is a degree $3$ polynomial if $h= \overline{f}$, $f$ is a cusp form; $P(\log x)$ is a degree $4$ polynomial if $h= E_r$, with $r\neq 0$; $P(\log x)$ is a degree $6$ polynomial if $h= E_0$. The constant in front of the leading term of $P(\log x)$ is given by \begin{equation*}
 \begin{cases}
\frac{L(1,f\otimes h)}{2\zeta (2)}, & \text{if $h\neq \overline{f}$, atleast one of $f,h$ is cuspidal,}\\
\frac{L(1,\mathrm{Ad}^2 f)}{\zeta (2)}, & \text{if $h= \overline{f}$, $f$ is cuspidal,}\\
\frac{|\zeta (1+2ir)|^2}{\zeta (2)}, & \text{if $h= E_r$, $r\neq 0$,}\\
\frac{1}{3\zeta (2)}, & \text{if $h= E_0$.}
\end{cases}   
\end{equation*} 
\end{theorem}

\noindent
As an immediate corollary we get the following corollary:

    
\begin{corollary}
    If we take $\Pi = f \boxplus h$, a $\mathrm{GL}(4)$ form, then for each $k$, large enough, there exists $g \in H_{k}$ with $K\leq k \leq 2K$ (see \cite{BXS}) such that,
    $$L(\frac{1}{2}, \Pi \otimes g)\neq 0 .$$
\end{corollary}

\begin{remark}
    For the case where $h=\overline{f}$, Blomer and Harcos (see \cite[Theorem 1]{BH2}) proved an asymptotic formula after applying an additional averaging over the critical line. Note that in that work, they considered an average of $| L(1 / 2 + it , f\otimes g)|^2$ over $t$ and $g$. In this article our main aim is to remove the average over $t$ and getting a asymptotic result of same strength.
\end{remark}

\begin{remark}
    The main ingredients of this paper is to use the approximate functional equation and the $\mathrm{GL}(2)$ Voronoi summation formula to reduce the original problem to a $\mathrm{GL}(2) \times \mathrm{GL}(2)$ shifted convolution sum problem and then appeal to some available bounds for those shifted sums.
\end{remark}

\begin{theorem}\label{shifted}
    Let $f,g$ be two automorphic forms for the full modular group with normalized Fourier coefficients $\lambda_{f}(n), \lambda_{g}(n)$ respectively. Let $0<h\ll N$ and $W$ be smooth nice weight function supported on $[1/2,3]$ with $W(x)=1$ for $x \in [1,2]$ and $W^{j}(x)\ll_{j}1$. If atleast one of $f$ or $g$ is a cusp form, we have
    $$\sum_{n=1}^{\infty}\lambda_{f}(n)\lambda_{g}(n+h)W\left( \frac{n}{N}\right)\ll_{\epsilon}N^{3/4}.$$
\end{theorem}

\begin{remark}
    This is a weaker bound than the square root bound got in \cite{BL} and \cite{BH}. But in those papers, they considered only cusp forms. But in our article, we are also considering non-cusp forms also.
\end{remark}

\begin{remark}
    If both of $f,g$ are not cusp forms then we will have a main term in the Theorem \ref{shifted}, coming from the Voronoi summation formula.
\end{remark}


\subsection*{Notation}
In this article, by $A \ll B$ we mean $|A| \leq C |B|$ for some absolute constant $C > 0$, depending only on $g,\epsilon$ and notation `$X\asymp Y$' will mean that $Y p^{-\epsilon}\leq X \leq Y p^{\epsilon}$. 

\subsection*{Acknowledgement}
The author thanks Gergely Harcos, Ritabrata Munshi and Rizwanur Khan for their helpful comments and suggestions. The author is also grateful to the Alfr\'ed R\'enyi Institute of Mathematics, Budapest, for the excellent research environment, as part of the work was done when the author was a young research scholar at the Alfr\'ed R\'enyi Institute of Mathematics. Also the author thanks to EPFL for the excellent research environment, and the author is supported by the project ``10003145 - Equidistribution in Number Theory".

\section{Sketch of the proof} 

\noindent
Here we will consider the case when $f,h$ are holomorphic cusp forms as the proof will be same for the Maa\ss \; cusp forms. Our sum is roughly

$$S:=\sum_{k \sim K}\sum_{g \in H_{k}} \omega_{g}^{-1}L(\frac{1}{2},f \otimes g)L(\frac{1}{2},g \otimes h).$$

\noindent
Using the Convexity result at this stage we will get the trivial bound, so at this stage we have a loss of size $K\times K = K^2$ for the $S$-sum.

\noindent
 After using the approximate functional equation, roughly our sum looks like

$$S=\sum_{k \sim K}\sum_{g \in H_{k}} \omega_{g}^{-1}\sum_{m \sim K^2}\lambda_{f}(m)\lambda_{g}(m)\sum_{n \sim K^2}\lambda_{h}(n)\lambda_{g}(n).$$

\noindent
Now we consider the sum
$$\sum_{g \in H_{k}} \omega_{g}^{-1}\lambda_{g}(m)\lambda_{g}(n).$$

\noindent
Applying the Petersson trace formula the above sum becomes

$$\delta_{n,m} + \sum_{c \ll K^2} i^k S(n,m;c) J_{k-1}\left( \frac{4\pi \sqrt{nm}}{c}\right),$$

\noindent
where $\delta_{n,m}$ is the Kronecker delta symbol. Hence the $S$ sum becomes

\begin{equation*}
    \begin{split}
        S \approx &\sum_{k \sim K}\sum_{n \sim K^2}\lambda_{f}(n)\lambda_{h}(n) \\
        &+ \sum_{k \sim K}\sum_{m \sim K^2}\lambda_{f}(m)\lambda_{g}(m)\sum_{n \sim K^2}\lambda_{h}(n)\lambda_{g}(n)\sum_{c \ll K^2} i^k S(n,m;c) J_{k-1}\left( \frac{4\pi \sqrt{nm}}{c}\right).
    \end{split}
\end{equation*}

\noindent
Now the first term (which is the diagonal term, getting after $n=m$), we note that the smooth sum $$\mathop{\sum}_{n \sim K^2}\lambda_{f}(n)\lambda_{h}(n) $$ is very small so the diagonal term is of size $K$. 

\noindent
For the off-diagonal term, i.e., for $n\neq m$, we have,

$$\sum_{k \sim K}\sum_{m \sim K^2}\lambda_{f}(m)\lambda_{g}(m)\sum_{n \sim K^2}\lambda_{h}(n)\lambda_{g}(n)\sum_{c \ll K^2} i^k S(n,m;c) J_{k-1}\left( \frac{4\pi \sqrt{nm}}{c}\right),$$

\noindent
  From the sum over $k$, we have $c\times K^2 \ll \sqrt{mn}\asymp K^2$, i.e., $c \sim K^\epsilon$. Here after applying the Petersson trace formula we save $\frac{\sqrt{K^2}}{\sqrt{c}} =K$ for the off-diagonal sum. So we need to save $\frac{K^2}{K}=K$ and a little more. At this stage the sum will look like

$$\mathop{\sum\sum}_{n,m\sim K^2}\lambda_{f}(n)\lambda_{h}(m)e(2\sqrt{mn}).$$

\noindent
Now first consider the $n$-sum which is

$$\sum_{n \sim K^2}\lambda_{f}(n)e(2\sqrt{mn}).$$

\noindent
Using the $\mathrm{GL}(2)$ Voronoi summation formula the above sum reduces to

$$\sum_{n \sim K^2}\lambda_{f}(n)e(2\sqrt{mn})\rightsquigarrow K^2\sum_{n \sim K^2}\lambda_{f}(n)\int e(2\sqrt{K^2 my}) J_{K-1}(4\pi \sqrt{nK^2 y}).$$

\noindent
Using the decomposition of the Bessel function the above sum reduces to

$$\sum_{n \sim K^2}\lambda_{f}(n)e(2\sqrt{mn})\rightsquigarrow K\sum_{\substack{n \sim K^2 \\ |n-m| \ll K^\epsilon}}\lambda_{f}(n).$$

\noindent
Hence we save $K$ here. So at this stage we save $K\times K =K^2$. Hence we are at the boundary and any savings will work. Now taking $n-m = d$ we note that we have to deal with the sum

$$\sum_{d\sim K^\epsilon}\sum_{n\sim K^2}\lambda_{f}(n)\lambda_{h}(n+d).$$

\noindent
Using the cancellation for the shifted convolution (see the Theorem \ref{shifted}) we will save $\sqrt{K}$. Hence the off-diagonal sum will be of size $\sqrt{K}$ and we will get our result.

\section{Preliminaries}
\subsection{Voronoi summation formula}
Here first we will recall the generalize Voronoi summation formula (see \cite{BH2}) for the Hecke eigenvalues $\lambda_{f}(n)$  twisted by a finite order additive character. 

\begin{proposition}\label{Vor}
\noindent
    Let $a$ and $c$ be coprime positive integers, and let $F:(0,\infty)\rightarrow \mathbb{C}$ be a smooth function of compact support. Then 
    $$c\sum_{n=1}^{\infty}\lambda_{f}(n)e\left( n\frac{a}{c}\right)F(n)=\sum_{n=1}^{\infty}\lambda_{f}(n) \sum_{\pm}e\left( \mp n \frac{\overline{a}}{c}\right)\int_{0}^{\infty}F(x)J_{f}\left( \frac{4\pi \sqrt{nx}}{c}\right) \; \mathrm{d}x,$$
    \noindent
    where
    $$J_{f}^{+}(x):=2\pi i^{k_{f}} J_{k_{f}-1}(x), \; J_{f}^{-}(x)=0,$$
    \noindent
    if $f$ is a holomorphic cusp form of level $1$ and weight $k_{f}$;
    $$J_{f}^{+}(x):= \frac{-\pi}{\cosh{(\pi r)}}\left( Y_{2ir}(x)+ Y_{-2ir}(x)\right) , \; J_{f}^{-}(x):=\epsilon_{f}4\cosh{(\pi r)}K_{2ir}(x),$$
    \noindent
    if $f$ is a Maa\ss \; cusp form of level $1$, weight $0$, Laplacian eigenvalue $\frac{1}{4}+r^2$ and sign $\epsilon_{f}\in \{ \pm 1\}$. For $f= E_{r}$ $(r \in \mathbb{R})$ the same formula holds with $J_{f}^{\pm}$ as in the Maa\ss \;  case (with $\epsilon_f =1$), except that on the right-hand side the following polar term has to be added:
$$\sum_{\pm}\zeta (1\pm 2ir) \int_{0}^{\infty}\left( \frac{x}{c^2}\right)^{\pm ir} F(x) \; \mathrm{d}x, \; \textit{ for } \; r\neq 0,$$
$$\int_{0}^{\infty}\left(\log \left(  \frac{x}{c^2}\right) +2\gamma \right) F(x) \; \mathrm{d}x, \; \textit{ for } \; r= 0.$$
\end{proposition}

\subsubsection{Bounds for the Bessel functions}\label{bessel} For $x\gg 1$, we have
$$J_{k-1}(2\pi x)= W_{k}(x)\frac{e(x)}{\sqrt{x}}+\overline{W}_{k}(x)\frac{e(-x)}{\sqrt{x}},$$

\noindent
for some function $W_{k}$ satisfying $x^j (\partial^j W_{k})(x)\ll_{j,k}1$.

\subsubsection{The Ramanujan bound on average} Here we note the following Ramanujan bound on average (see \cite{IW}):
\begin{equation*}
    \sum_{n\ll N}|\lambda_{f}(n)|^2 \ll  N^{1+\epsilon}.
\end{equation*}

\subsection{$\mathrm{GL}(2)$ holomorphic cusp forms} We will use the following well-known (\cite{Iw}) Petersson trace formula
\begin{proposition}\label{Pet}
   $$ \sum_{g \in H_{k}}\omega_{g}^{-1}\lambda_{g}(m)\lambda_{g}(n)= \delta_{m,n}+ 2\pi i^{-k}\sum_{c=1}^{\infty}\frac{S(m,n;c)}{c}J_{k-1}\left( \frac{4\pi \sqrt{mn}}{c}\right),$$
   where $\omega_{g}^{-1}= \frac{(4\pi)^{k-1}}{\Gamma (k-1)}\lVert g \rVert^{2}$, $\delta_{m,n}$ equals to $1$ if $m=n$ and $0$ otherwise, $S(m,n;c)$ is the Kloosterman sum defined below, and $J_{k-1}$ is the $J$-Bessel function.
\end{proposition}

\noindent
The Kloosterman sum is defined as :
$$S(m,n;c) =\sum_{d\overline{d}\equiv 1\bmod c}e\left( \frac{md + n \overline{d}}{c}\right) .$$
\noindent
Weil proved that
$$\lvert S(m,n;c)\rvert \leq (m,n,c)^{1/2}c^{1/2}\tau (c) ,$$
where $\tau$ is the divisor function. Then clearly one can see that

\begin{equation*}
    \sum_{n \leq x}|S(m,n;c)| \ll x c^{1/2 +\epsilon}.
\end{equation*}

\subsection{Rankin-Selberg $L$-function} For $f \in H_{\ell}$, $g \in H_{k}$, $h\in H_{m}$, we will consider the Rankin-Selberg $L$-functions $L(s, f\otimes g)$ and $L(s, g\otimes h)$ where
$$L(s, f\otimes g) = \zeta (2s)\sum_{m=1}^{\infty}\frac{\lambda_{f}(m)\lambda_{g}(m)}{m^s} = \sum_{m=1}^{\infty}\frac{\lambda_{f\otimes g}(m)}{m^s},$$
where 
$$\lambda_{f\otimes g}(m)= \sum_{a b^2 = m}\lambda_{f}(a)\lambda_{g}(a) .$$
This becomes an entire function and satisfies the following functional equation
$$\Lambda(s, f\otimes g)= \Lambda (1-s, f\otimes g)$$
where
$$\Lambda(s, f\otimes g)= \frac{1}{(4\pi^2)^s}\Gamma \left( s+ \frac{k+\ell}{2}-1\right)\Gamma \left(s+ \frac{k-\ell}{2} \right)L(s, f\otimes g).$$
Then we have the following Approximate Functional equation :
\begin{lemma}\label{afe}
    Let $G(u)= e^{u^2}$. Then we have
    \begin{equation}
        L(\frac{1}{2},f\otimes g)= 2 \sum_{m=1}^{\infty}\frac{\lambda_{f\otimes g}(m)}{m^{1/2}}U(m,k),
    \end{equation}
    where
    \begin{equation}\label{U}
        U(y,k)= \frac{1}{2\pi i}\int_{(3)}y^{-u}G(u)\frac{\gamma(\frac{1}{2}+u,k)}{\gamma(\frac{1}{2},k)}\frac{\mathrm{d}u}{u},
    \end{equation}
    and
    \begin{equation}
        \gamma(s,k)= (2\pi)^{-2s}\Gamma\left( s+\frac{k+\ell}{2}-1\right)\Gamma\left( s+\frac{k-\ell}{2}\right) .
    \end{equation}
\end{lemma}
\begin{proof}
    For the proof, see \cite[Theorem $5.3$]{IK}.
\end{proof}
\begin{lemma}\label{UV}
    $(1)$ The derivatives of $U(y,k)$ with respect to $y$ satisfy
    $$y^a \frac{\partial^a}{\partial y^a}U(y,k)\ll \left(1 + \frac{y}{k^2} \right)^{-A} ,$$
    where the implied constants depends only on $a$ and $A$.
    $(2)$ For $1 \leq y \ll k^{2+\epsilon}$, as $k \mapsto \infty$, we have
    $$U(y,k)= U_{1}\left( \frac{y}{k^2}\right)+ O\left( \frac{1}{k}U_{1}^{*}\left(\frac{y}{k^2} \right)\right) ,$$
    where
    $$U_{1}(x)=\frac{1}{2\pi i}\int_{(3)}x^{-u}G_{2}(u) \frac{\mathrm{d}u}{u}\, \, , \, \, U_{1}^{*}(x)=\frac{1}{2\pi i}\int_{(3)}x^{-u}G_{3}(u)\frac{\mathrm{d}u}{u},$$
    and $G_{2}(u)$, $G_{3}(u)$ are holomorphic functions with exponential decay as $\lvert \mathcal{I}u\rvert \mapsto \infty$.
\end{lemma}
\begin{proof}
\noindent
    $(1)$ See \cite[Proposition $5.4$]{IK}.

    $(2)$ By Stirling's formula we have,
    \begin{align*}
        G(u)\frac{\gamma \left(\frac{1}{2}+u,k \right)}{\gamma \left(\frac{1}{2},k \right)} & = G(u)\left( \frac{k^2}{16 \pi^2}\right)^{u}\left(1+ O\left( \frac{p(u)}{k}\right)\right) \\
        & = k^{2u} G_{2}(u) + O\left( \frac{1}{k}k^{2u}G_{3}(u)\right),
    \end{align*}
    where $p(u)$ is a polynomial in $u$ with $G_{2}(u)= \frac{G(u)}{(16\pi^2 )^u }$ and $G_{3}(u) = \frac{G(u) p(u)}{(16 \pi^2 )^u}$. Hence inserting these into \eqref{U}, we get our result.
\end{proof}
\section{Proof of theorem \ref{main}}

\noindent
We have, by the approximate functional equation and the Petersson trace formula, we have

\begin{align*}
    &\sum_{g \in H_{k}}\omega_{g}^{-1}L(\frac{1}{2}, f\otimes g)L(\frac{1}{2}, g\otimes h) \\
    & = 4 \sum_{g \in H_{k}}\omega_{g}^{-1}\sum_{m=1}^{\infty}\frac{\lambda_{f\otimes g}(m)}{m^{1/2}}U(m,k)\sum_{n=1}^{\infty}\frac{\lambda_{g\otimes h}(n)}{n^{1/2}}V(n,k) \\
    & = 4 \sum_{g \in H_{k}}\omega_{g}^{-1}\sum_{a_{1}=1}^{\infty}\sum_{b_{1}=1}^{\infty}\frac{\lambda_{f}(a_1 )\lambda_{g}(a_{1})}{(a_{1}b_{1}^2)^{1/2}}U(a_1 b_{1}^{2},k)\sum_{a_{2}=1}^{\infty}\sum_{b_{2}=1}^{\infty}\frac{\lambda_{g}(a_2 )\lambda_{h}(a_{2})}{(a_{2}b_{2}^2)^{1/2}}V(a_2 b_{2}^{2},k) \\
    & = 4 \sum_{a_{1}=1}^{\infty}\sum_{b_{1}=1}^{\infty}\sum_{b_{2}=1}^{\infty}\frac{\lambda_{f}(a_{1})\lambda_{h}(a_{1})}{a_{1}b_{1}b_{2}}U(a_{1}b_{1}^{2},k)V(a_{1}b_{2}^{2},k) + \\
    & 8\pi \sum_{a_{1}=1}^{\infty}\sum_{a_{2}=1}^{\infty}\sum_{b_{1}=1}^{\infty}\sum_{b_{2}=1}^{\infty}\frac{\lambda_{f}(a_{1})\lambda_{h}(a_{2})}{(a_{1}a_{2}b_{1}^{2}b_{2}^{2})^{1/2}}U(a_{1}b_{1}^{2},k)V(a_{2}b_{2}^{2},k)\sum_{c=1}^{\infty}i^{k}\frac{S(a_{1},a_{2};c)}{c}J_{k-1}\left(\frac{4\pi \sqrt{a_{1}a_{2}}}{c} \right) \\
    & := 4 D_{k}+ 8\pi E_{k}.
\end{align*}
    
\subsection{Diagonal terms}\label{diagonal} First we need to deal with the diagonal terms (i.e., when $a_{1}=a_{2}$) which contributes in the main term :
$$D : = 4 \sum_{k \equiv 0 \bmod 2}u\left(\frac{k-1}{K} \right) D_{k}.$$

\begin{lemma}\label{poly}
    We have
    \begin{equation}\label{Dk}
       D_{k}= P(\log k) + \mathrm{O}(k^{-1}) .
    \end{equation}
   and
    \begin{equation}\label{D}
        \begin{split}
            D=2KP(\log K) + \mathrm{O}(K^\epsilon ) .
        \end{split}
    \end{equation}
    \noindent
    where $P(\log x)$ is a polynomial in $\log x$, depending on $f,h$ only. $P(\log x)$ is a degree $2$ polynomial if $f\neq \overline{h}$, atleast one of $f,h$ is a cusp form; $P(\log x)$ is a degree $3$ polynomial if $f= \overline{h}$, $f$ is a cusp form; $P(\log x)$ is a degree $4$ polynomial if $h= E_r$, with $r\neq 0$; $P(\log x)$ is a degree $6$ polynomial if $h= E_0$. The constant in front of the leading term of $P(\log x)$ is given by \eqref{constant}. 
    
\end{lemma}
\begin{proof}
\begin{equation*}
    \begin{split}
       D_{k}&= \sum_{a_{1}=1}^{\infty}\sum_{b_{1}=1}^{\infty}\sum_{b_{2}=1}^{\infty}\frac{\lambda_{f}(a_{1})\lambda_{h}(a_{1})}{a_{1}b_{1}b_{2}}U(a_{1}b_{1}^{2},k)V(a_{1}b_{2}^{2},k) \\
       &= \frac{1}{(2\pi i)^{2}}\iint_{(3)}\\
&\times\sum_{a_{1}=1}^{\infty}\frac{\lambda_{f}(a_{1})\lambda_{h}(a_{1})}{a_{1}^{1+v+w}}\sum_{b_{1}=1}^{\infty}\frac{1}{b_{1}^{1+2v}}\sum_{b_{2}=1}^{\infty}\frac{1}{b_{2}^{1+2w}}\frac{\gamma\left(\frac{1}{2}+v,k \right)}{\gamma\left(\frac{1}{2},k \right)}\frac{\gamma\left(\frac{1}{2}+w,k \right)}{\gamma\left(\frac{1}{2},k \right)}\frac{G(v)}{v}\frac{G(w)}{w} \, \mathrm{d}v \, \mathrm{d}w \\
       &= \frac{1}{(2\pi i)^{2}}\iint_{(3)}\\
       &\times\frac{L(1+v+w, f\otimes h)}{\zeta (2+2v+2w)}\zeta (1+2v)\zeta (1+2w)\frac{\gamma\left(\frac{1}{2}+v,k \right)}{\gamma\left(\frac{1}{2},k \right)}\frac{\gamma\left(\frac{1}{2}+w,k \right)}{\gamma\left(\frac{1}{2},k \right)}\frac{G(v)}{v}\frac{G(w)}{w} \, \mathrm{d}v \, \mathrm{d}w .
    \end{split}
\end{equation*}

First moving the line of integration to $\mathcal{R}(v)=-2$, we note that we get a double pole at $v=0$. Then by the Residue theorem we get that
\begin{equation*}
    \begin{split}
        D_k &= \frac{1}{2\pi i}\int_{(3)}\frac{L^{\prime}(1+w,f\otimes h)}{2\zeta (2+2w)}\zeta (1+2w)\frac{\gamma\left(\frac{1}{2}+w,k\right)}{\gamma\left(\frac{1}{2},k\right)}\frac{G(w)}{w} \; \mathrm{d}w \\
        & - \frac{1}{2\pi i}\int_{(3)}\frac{L(1+w,f\otimes h)}{2(\zeta (2+2w))^2}\zeta (1+2w)\frac{\gamma\left(\frac{1}{2}+w,k\right)}{\gamma\left(\frac{1}{2},k\right)}\frac{G(w)}{w} \; \mathrm{d}w
     \\
     & + \frac{1}{2\pi i}\int_{(3)}\frac{\gamma L(1+w,f\otimes h)}{\zeta (2+2w)}\zeta (1+2w)\frac{\gamma\left(\frac{1}{2}+w,k\right)}{\gamma\left(\frac{1}{2},k\right)}\frac{G(w)}{w} \; \mathrm{d}w \\
     &+ \frac{1}{2\pi i}\int_{(3)}\frac{L(1+w,f\otimes h)}{2\zeta (2+2w)}\zeta (1+2w)\frac{\gamma^{\prime}\left(\frac{1}{2}+w,k\right)}{\gamma\left(\frac{1}{2},k\right)}\frac{G(w)}{w} \; \mathrm{d}w \\
     &+ \frac{1}{(2\pi i)^{2}}\int_{(-2)}\int_{(3)}\\
       &\times\frac{L(1+v+w, f\otimes h)}{\zeta (2+2v+2w)}\zeta (1+2v)\zeta (1+2w)\frac{\gamma\left(\frac{1}{2}+v,k \right)}{\gamma\left(\frac{1}{2},k \right)}\frac{\gamma\left(\frac{1}{2}+w,k \right)}{\gamma\left(\frac{1}{2},k \right)}\frac{G(v)}{v}\frac{G(w)}{w} \, \mathrm{d}v \, \mathrm{d}w ,
    \end{split}
\end{equation*}

\noindent
where $\gamma$ is the Euler-Mascheroni constant. Now we will move the line of integration to $\mathcal{R}(w)= -\frac{1}{2}$ so that in that process we will get a pole at $w=0$ of order $e$ which depends on $h$: if $h$ is a cusp form and $h \neq \overline{f}$, then $e = 2$; if $h=\overline{f}$, then $e=3$; if $h=E_r$ with $r\neq 0$, then $e=4$; while for $h=E_{0}$ we have $e=6$.

Then we can note that the residue of the pole is given by the linear combination of $L^{(j)}(1,f\otimes h), \frac{\gamma^{(j)}\left(\frac{1}{2},k \right)}{\gamma\left(\frac{1}{2},k \right)}$ and other constants, not depending on $f,h$ and also integrals on the lines $\mathcal{R}(v)=-2$, $\mathcal{R}(w)=-\frac{1}{2}$ where $j \leq e$. By Stirling's formula, we have
$$\frac{\gamma^{\prime}\left(\frac{1}{2} ,k\right)}{\gamma\left( \frac{1}{2},k\right)} = 2 \log k -4 \log 2 - 2 \log \pi + O(k^{-1}),$$
$$\frac{\gamma^{(j)}\left(\frac{1}{2} ,k\right)}{\gamma\left( \frac{1}{2},k\right)} = (\log k^2 )^j + \mathrm{O}((\log k^2 )^{j-1}),$$
$$\frac{\gamma\left(\frac{1}{2}+v, k \right)}{\gamma\left(\frac{1}{2},k\right)} \ll k^{-2} ,$$
for $\mathcal{R}(v)=-2$, and
$$\frac{\gamma\left(\frac{1}{2}+w, k \right)}{\gamma\left(\frac{1}{2},k\right)} \ll k^{-1} ,$$
for $\mathcal{R}(w)=-\frac{1}{2}$. 

Hence after shifting the line of integration to $\mathcal{R}(w)= -\frac{1}{2}$, the leading term will come from the term

$$\frac{1}{2\pi i}\int_{(3)}\frac{L(1+w,f\otimes h)}{2\zeta (2+2w)}\zeta (1+2w)\frac{\gamma^{\prime}\left(\frac{1}{2}+w,k\right)}{\gamma\left(\frac{1}{2},k\right)}\frac{G(w)}{w} \; \mathrm{d}w .$$

\noindent
The constant in front of the leading term equals
\begin{equation}\label{constant}
 \frac{\displaystyle\lim_{w\rightarrow 0}\left(2^{e-2} w^{e-2}L(1+w, f\otimes h)\right)}{4\zeta (2) \; (e-2)!}=\begin{cases}
\frac{L(1,f\otimes h)}{4\zeta (2)}, & \text{if $h\neq \overline{f}$, atleast one of $f,h$ is cuspidal,}\\
\frac{L(1,\mathrm{Ad}^2 f)}{2\zeta (2)}, & \text{if $h= \overline{f}$, $f$ is cuspidal,}\\
\frac{|\zeta (1+2ir)|^2}{2\zeta (2)}, & \text{if $h= E_r$, $r\neq 0$,}\\
\frac{1}{6\zeta (2)}, & \text{if $h= E_0$.}
\end{cases}   
\end{equation}

Hence using the reasoning given above and the estimates done above, we note that , we get \eqref{Dk}.

\noindent
By the Poisson summation formula, we have,
\begin{equation}\label{Poi}
    4\sum_{k \equiv 0 \bmod 2}u\left(\frac{k-1}{K}\right)= 2K \int_{0}^{\infty}u(\xi ) \, \mathrm{d}\xi \, \, + O(K^{-A}),
\end{equation}
for any $A > 0$. So the main term will be of size $K$. Now by writing $\log k = \log \left(\frac{k-1}{K}K\right) + O(k^{-1})$ and using \eqref{Poi}, we get \eqref{D}.

\end{proof}

\begin{remark}
    We note that $L(1,f\otimes h)$ is never zero. In general, Shahidi in \cite{SH} proved that a $\mathrm{GL}(m)\times \mathrm{GL}(n)$ Rankin-Selberg $L$-function never vanishes on the line $\mathrm{Re}(s)=1$. Moreover, Thorner and Harcos in \cite{HT} proved that such an $L$-function does not vanish for

$$\mathrm{Re}(s) > 1-\frac{c}{(1+|\mathrm{Im}(s)|)^{\varepsilon}},$$

\noindent
where $\varepsilon>0$ is arbitrary, and $c$ is a constant depending on $\varepsilon$ and the $L$-function at hand.
\end{remark}

\subsection{Off-diagonal terms}
Now we want to show that the off-diagonal terms (i.e., when $a_{1}\neq a_{2}$) only contribute to the error term. Let
$$E := \sum_{2 |k}u\left( \frac{k-1}{K}\right)E_{k} $$
$$=\sum_{a_{1}=1}^{\infty}\sum_{a_{2}=1}^{\infty}\sum_{b_{1}=1}^{\infty}\sum_{b_{2}=1}^{\infty}\frac{\lambda_{f}(a_{1})\lambda_{h}(a_{2})}{(a_{1}a_{2}b_{1}^{2}b_{2}^{2})^{1/2}}\sum_{c=1}^{\infty}\frac{S(a_{1},a_{2};c)}{c}$$
$$\times \sum_{2 |k}i^{k}u\left( \frac{k-1}{K}\right)U(a_{1}b_{1}^{2},k)V(a_{2}b_{2}^{2},k)J_{k-1}\left(\frac{4\pi \sqrt{a_{1}a_{2}}}{c} \right)$$

\begin{lemma}\label{E}
    We have
    $$E \ll \frac{1}{K}.$$
\end{lemma}
 \begin{proof}
     At first we take a smooth dyadic subdivision, so that we need to estimate
     $$E_{N,M}:=\sum_{a_{1}=1}^{\infty}\sum_{a_{2}=1}^{\infty}\sum_{b_{1}=1}^{\infty}\sum_{b_{2}=1}^{\infty}w\left( \frac{a_{1}b_{1}^{2}}{N}\right)H\left( \frac{a_{2}b_{2}^{2}}{M}\right)\frac{\lambda_{f}(a_{1})\lambda_{h}(a_{2})}{(a_{1}a_{2}b_{1}^{2}b_{2}^{2})^{1/2}}\sum_{c=1}^{\infty}\frac{S(a_{1},a_{2};c)}{c}$$
     $$\times \sum_{2 |k}i^{k}u\left( \frac{k-1}{K}\right)U(a_{1}b_{1}^{2},k)V(a_{2}b_{2}^{2},k)J_{k-1}\left(\frac{4\pi \sqrt{a_{1}a_{2}}}{c} \right)$$
     for $N\ll K^{2+\epsilon}$ and $M\ll K^{2+\epsilon}$, where $w,h$ are fixed smooth functions with support contained in $[1/2,3]$ and $w(x)=h(x)=1$ whenever $x\in [1,2]$.

\noindent
     Let
     $$R_{1}:=\sum_{a_{1}=1}^{\infty}\sum_{a_{2}=1}^{\infty}\sum_{b_{1}=1}^{\infty}\sum_{b_{2}=1}^{\infty}w\left( \frac{a_{1}b_{1}^{2}}{N}\right)H\left( \frac{a_{2}b_{2}^{2}}{M}\right)\frac{\lambda_{f}(a_{1})\lambda_{h}(a_{2})}{(a_{1}a_{2}b_{1}^{2}b_{2}^{2})^{1/2}}\sum_{c\geq 16\pi \frac{\sqrt{NM}}{b_{1}b_{2}}}\frac{S(a_{1},a_{2};c)}{c}$$
     $$\times \sum_{2 |k}i^{k}u\left( \frac{k-1}{K}\right)U(a_{1}b_{1}^{2},k)V(a_{2}b_{2}^{2},k)J_{k-1}\left(\frac{4\pi \sqrt{a_{1}a_{2}}}{c} \right) ,$$

\noindent
and

$$R_{2}:=\sum_{a_{1}=1}^{\infty}\sum_{a_{2}=1}^{\infty}\sum_{b_{1}=1}^{\infty}\sum_{b_{2}=1}^{\infty}w\left( \frac{a_{1}b_{1}^{2}}{N}\right)H\left( \frac{a_{2}b_{2}^{2}}{M}\right)\frac{\lambda_{f}(a_{1})\lambda_{h}(a_{2})}{(a_{1}a_{2}b_{1}^{2}b_{2}^{2})^{1/2}}\sum_{c < 16\pi \frac{\sqrt{NM}}{b_{1}b_{2}}}\frac{S(a_{1},a_{2};c)}{c}$$
     $$\times \sum_{2 |k}i^{k}u\left( \frac{k-1}{K}\right)U(a_{1}b_{1}^{2},k)V(a_{2}b_{2}^{2},k)J_{k-1}\left(\frac{4\pi \sqrt{a_{1}a_{2}}}{c} \right) .$$

So we have $E_{N,M}= R_{1}+ R_{2}$. Now our aim is to show that $R_{1}$ is negligible and $R_{2}\ll K^{1/2 +\epsilon}$ which we will show in the following sections.
     
 \end{proof}

\section{Estimation of $R_{1}$}
\noindent
 For $c \geq 16\pi \frac{\sqrt{NM}}{b_{1}b_{2}}$, $\frac{N}{b_{1}^{2}}\leq a_{1} \leq \frac{2N}{b_{1}^{2}}$ and $\frac{M}{b_{2}^{2}}\leq a_{2} \leq \frac{2M}{b_{2}^{2}}$, we have
 $$x =\frac{4\pi \sqrt{a_{1} a_{2}}}{c}\leq \frac{1}{4}\left( \frac{a_{1}a_{2}}{NM}\right)^{1/2} < e^{-2} .$$
So we have $J_{k-1}(x)\ll x^{k-1}$ for $0 < x<1$. Now by the Weil's bound for Kloosterman sums we have
$$S(a_{1},a_{2};c)\ll c^{1/2 + \epsilon}(a_{1},a_{2},c)^{1/2} ,$$
and by the Ramanujan bound on average,

$$\sum_{n \leq x}|\lambda_{f} (n)|^2 \ll x .$$

\noindent
Using these informations, we have,

\begin{equation}\label{R1}
\begin{split}
    R_{1} & \ll \sum_{a_{1}=1}^{\infty}\sum_{a_{2}=1}^{\infty}\sum_{b_{1}=1}^{\infty}\sum_{b_{2}=1}^{\infty}\Big| w\left( \frac{a_{1}b_{1}^{2}}{N}\right)\Big| \, \Big|H\left( \frac{a_{2}b_{2}^{2}}{M}\right)\Big| \, \frac{|\lambda_{f}(a_{1})|}{(a_{1}b_{1}^{2})^{1/2}}\frac{|\lambda_{h}(a_{2})|}{(a_{2}b_{2}^{2})^{1/2}} \\
    & \times \sum_{c\geq 16\pi \frac{\sqrt{NM}}{b_{1}b_{2}}}c^{- 1/2 +\epsilon}(a_{1},a_{2},c)^{1/2}
     \times \sum_{2 |k}\Big|u\left( \frac{k-1}{K}\right)\Big| U(a_{1}b_{1}^{2},k)V(a_{2}b_{2}^{2},k)J_{k-1}\left(\frac{4\pi \sqrt{a_{1}a_{2}}}{c} \right) \\
     & \ll \sum_{N \leq a_{1}b_{1}^{2}\leq 2N}\frac{|\lambda_{f}(a_{1})|}{(a_{1}b_{1}^{2})^{1/2}}\sum_{N \leq a_{1}b_{1}^{2}\leq 2N}\frac{|\lambda_{h}(a_{2})|}{(a_{2}b_{2}^{2})^{1/2}}\sum_{c \geq 16\pi \frac{\sqrt{NM}}{b_{1}b_{2}}}\left( \frac{4\pi \sqrt{a_{1}a_{2}}}{c}\right)^{2} \sum_{2 | k, k \sim K} e^{-(k-3)} \\
     & \ll K^{-A},
\end{split}
\end{equation}
\noindent
for any $A > 0$.

\section{Estimation of $R_{2}$}

\begin{lemma}
    For a fixed function $g(\xi) \in \mathcal{C}_{c}^{\infty}(0,\infty)$, we have
    $$\sum_{2 | k}i^{k}g(k-1)J_{k-1}(x)=  F(x),$$
    where
    $$F(x)=\int_{-\infty}^{\infty}\hat{g}(t)\sin{(x\cos{(2\pi t)})} \, \mathrm{d}t ,$$
    and
    $$\hat{g}(t)= \int_{-\infty}^{\infty}g(\xi) e(t\xi ) \, \mathrm{d}\xi ,$$
    the Fourier transform of $g$.
\end{lemma}
\begin{proof}
    See \cite[P. $85-87$]{Iw}.
\end{proof}

\noindent
Now here we note that

$$F(x)=\frac{1}{2i}\left(\int_{-\infty}^{\infty}\hat{g}(t)e\left(\frac{x}{2\pi}\cos{(2\pi t)}\right) \, \mathrm{d}t - \, \int_{-\infty}^{\infty}\hat{g}(t)e\left(-\frac{x}{2\pi}\cos{(2\pi t)}\right) \, \mathrm{d}t \right),$$

\noindent
so without loss of generality we write

\begin{equation*}
    F(x)=\int_{-\infty}^{\infty}\hat{g}(t)e\left(\pm\frac{x}{2\pi}\cos{(2\pi t)}\right) \, \mathrm{d}t .
\end{equation*}

\noindent
Now let
$$g(\xi )= u\left( \frac{\xi}{K}\right)U(a_{1}b_{1}^{2},\xi +1)V(a_{2}b_{2}^{2},\xi +1),$$

\noindent
and 

$$g_{0}(\xi )= u\left( \xi\right)U(a_{1}b_{1}^{2},(K\xi +1)^2 )V(a_{2}b_{2}^{2},(K\xi +1)^2 ).$$

\noindent
Then from the lemma \ref{UV} we need to consider the contribution only from

$$u\left( \xi\right)U_{1}\left( \frac{a_{1}b_{1}^{2}}{(K\xi +1)^2 }\right) V_{1}\left( \frac{a_{2}b_{2}^{2}}{(K\xi +1)^2}\right),$$

\noindent
so by abuse of notation we can write

\begin{equation}\label{g0}
  g_{0}(\xi )=  u\left( \xi\right)U_{1}\left( \frac{a_{1}b_{1}^{2}}{(K\xi +1)^2}\right) V_{1}\left( \frac{a_{2}b_{2}^{2}}{(K\xi +1)^2}\right).
\end{equation}

\noindent
Then by changing the variables we can see that $\hat{g}(t)= K \hat{g}_{0}(Kt)$ so that we have

\begin{equation}\label{F}
    F(x)=\int_{-\infty}^{\infty}\hat{g}_{0}(t)e\left(\pm\frac{x}{2\pi}\cos{\left(\frac{2\pi t}{K}\right)}\right) \, \mathrm{d}t .
\end{equation}

\subsection{Simplification of $R_{2}$}
 For the estimation of $R_{2}$, we need the Voronoi summation formula \ref{Vor}. Now consider the sum

\begin{equation*}
    \begin{split}
R_{2}&= \sum_{a_{1}=1}^{\infty}\sum_{a_{2}=1}^{\infty}\sum_{b_{1}=1}^{\infty}\sum_{b_{2}=1}^{\infty}w\left( \frac{a_{1}b_{1}^{2}}{N}\right)H\left( \frac{a_{2}b_{2}^{2}}{M}\right)\frac{\lambda_{f}(a_{1})\lambda_{h}(a_{2})}{(a_{1}a_{2}b_{1}^{2}b_{2}^{2})^{1/2}}U(a_{1}b_{1}^{2},k)V(a_{2}b_{2}^{2},k)\\
&\times \sum_{c < 16\pi \frac{\sqrt{NM}}{b_{1}b_{2}}}\frac{S(a_{1},a_{2};c)}{c}F\left( \frac{4\pi \sqrt{a_{1}a_{2}}}{c}\right) .
    \end{split}
\end{equation*}

\noindent
First consider the function given in \eqref{F}:

\begin{equation*}
    F(x)=\int_{-\infty}^{\infty}\hat{g}_{0}(t)e\left(\pm\frac{x}{2\pi}\cos{\left(\frac{2\pi t}{K}\right)}\right) \, \mathrm{d}t .
\end{equation*}

\noindent
Then here as the contribution of $|t| \geq K^\epsilon$ is negligible, we need to consider the case when $|t| \leq K^\epsilon$.  Then expanding $\cos{\left( \frac{2\pi t}{K}\right)}$ into the Taylor series and then using $\doublehat{g}_{0}(x)= g_{0}(-x)$, we get that

\begin{equation*}
    \begin{split}
        F(x)& = \int_{-\infty}^{\infty}\hat{g}_{0}(t)e\left(\frac{x}{2\pi} -\frac{\pi xt^2}{K^2}\right)\mathrm{d}t + \mathrm{O}\left(\int_{-\infty}^{\infty}|\hat{g}_{0}(t)|\frac{|x|}{K^4}|t|^4 \, \mathrm{d}t\right)\\
        &= e\left(\frac{x}{2\pi}\right)\int_{-\infty}^{\infty}\hat{g}_{0}(t)e\left( -\frac{\pi xt^2}{K^2}\right)\mathrm{d}t + \mathrm{O}\left(\frac{|x|}{K^4} \, \right),
    \end{split}
\end{equation*}

\noindent
as we have $$\widehat{g_{0}^{(a)}}(t)= (-2\pi i t)^a \, \widehat{g_{0}}(t).$$

\noindent
Then by the Weil bound for the Kloosterman sum and the Ramanujan bound on average, the above $\mathrm{O}$-term contributes at most

$$\frac{NM}{K^4} \ll K^{\epsilon},$$

\noindent
so without loss of generality, by abuse of notation we have

\begin{equation}\label{aF}
    \begin{split}
        F(x)& = e\left(\frac{x}{2\pi}\right)\int_{-\infty}^{\infty}\hat{g}_{0}(t)e\left( -\frac{\pi xt^2}{K^2}\right)\mathrm{d}t \\
        & = \frac{K}{\sqrt{x}}e\left(\frac{x}{2\pi}\right)\int_{-\infty}^{\infty}g_{0}(-t)e\left( -\frac{\pi K^2 t^2}{x}\right)\mathrm{d}t \\
        & = \frac{K}{\sqrt{x}}e\left(\frac{x}{2\pi}\right)\int_{-\infty}^{\infty}g_{0}(t)e\left( -\frac{\pi K^2 t^2}{x}\right)\mathrm{d}t ,
    \end{split}
\end{equation}

\noindent
where $g_{0}$ is given by \eqref{g0} and for the second line we have used the Parseval's identity. Now by repeated integration by parts, we note that the above integral becomes negligibly small unless $|x|\sim K^2$. Now for $|x|= \frac{4\pi \sqrt{a_{1}a_{2}}}{c}$, we see that $c\sim \frac{4\pi \sqrt{a_{1}a_{2}}}{K^2}\sim K^\epsilon $. Now as $|x|\sim K^2$ so without loss of generality we have
$$\int_{-\infty}^{\infty}g_{0}(t)e\left( -\frac{\pi K^2 t^2}{x}\right)\mathrm{d}t = A,$$

\noindent
for some constant $A>0$.Then the sum $R_{2}$ reduces to
\begin{equation*}
    \begin{split}
R_{2}&= \sum_{a_{1}=1}^{\infty}\sum_{a_{2}=1}^{\infty}\sum_{b_{1}=1}^{\infty}\sum_{b_{2}=1}^{\infty}w\left( \frac{a_{1}b_{1}^{2}}{N}\right)H\left( \frac{a_{2}b_{2}^{2}}{M}\right)\frac{\lambda_{f}(a_{1})\lambda_{h}(a_{2})}{(a_{1}a_{2}b_{1}^{2}b_{2}^{2})^{1/2}}\\
&\times \sum_{c \sim K^\epsilon}\frac{S(a_{1},a_{2};c)}{c} \frac{K\sqrt{c}}{(a_{1}a_{2})^{1/4}}e\left( \frac{2 \sqrt{a_{1}a_{2}}}{c}\right) \\
& = K\sum_{b_{1}=1}^{\infty}\sum_{b_{2}=1}^{\infty}\frac{1}{b_{1}b_{2}}\sum_{c \sim K^\epsilon}\frac{1}{\sqrt{c}} \, \, \sideset{}{^*}\sum_{\beta \bmod c}\, \sum_{a_{2}=1}^{\infty}\frac{\lambda_{h}(a_{2})}{a_{2}^{3/4}}e\left(\frac{a_{2}\beta}{c} \right)H\left( \frac{a_{2}b_{2}^{2}}{M}\right)\\
& \times \left(\sum_{a_{1}=1}^{\infty} \frac{\lambda_{f}(a_{1})}{a_{1}^{3/4}}e\left(\frac{a_{1}\overline{\beta}}{c} \right)e\left( \frac{2 \sqrt{a_{1}a_{2}}}{c}\right) w\left( \frac{a_{1}b_{1}^{2}}{N}\right)\right).
    \end{split}
\end{equation*}

\noindent
Now consider the last sum

$$\sum_{a_{1}=1}^{\infty} \frac{\lambda_{f}(a_{1})}{a_{1}^{3/4}}e\left(\frac{a_{1}\overline{\beta}}{c} \right)e\left( \frac{2 \sqrt{a_{1}a_{2}}}{c}\right) w\left( \frac{a_{1}b_{1}^{2}}{N}\right) .$$

Applying the Voronoi summation formula \ref{Vor} the above sum reduces to the following sum:

$$\frac{2\pi i^k}{c}\sum_{a_{1}=1}^{\infty} \lambda_{f}(a_{1}) e\left(\frac{- a_{1}\beta}{c} \right)\mathcal{I}_{1} ,$$

\noindent
where the integral $\mathcal{I}_{1}$ is given by

$$\mathcal{I}_{1}=\int_{0}^{\infty}e\left( \frac{2\sqrt{xa_{2}}}{c}\right)\frac{1}{x^{3/4}}w\left( \frac{x b_{1}^{2}}{N}\right) \, J_{k-1}\left( \frac{4\pi \sqrt{xa_{1}}}{c}\right) \, \, \mathrm{d}x.$$

\noindent
Changing the variables $\frac{x b_{1}^{2}}{N} \mapsto x$ we have

$$\mathcal{I}_{1}=\frac{N^{1/4}}{\sqrt{b_{1}}}\int_{0}^{\infty}e\left( \frac{2 \sqrt{Nxa_{2}}}{cb_{1}}\right)\frac{1}{x^{3/4}}w\left( x\right) \, J_{k-1}\left( \frac{4\pi \sqrt{xNa_{1}}}{cb_{1}}\right) \, \, \mathrm{d}x.$$

\noindent
Using the decomposition, 
 $$J_{k -1}(2\pi x)= \frac{W (x)}{\sqrt{x}}e(x)+ \frac{\bar{W }(x)}{\sqrt{x}}e(-x) ,$$

 \noindent
  where $W (x)$ is a nice function, the integral becomes

$$\mathcal{I}_{1}=\frac{\sqrt{c}}{ a_{1}^{1/4}}\int_{0}^{\infty}e\left( \frac{2 \sqrt{Nxa_{2}}}{cb_{1}}\right)w\left( x\right) \frac{W(x)}{x} \, e\left( \mp \frac{2 \sqrt{xNa_{1}}}{b_{1}c}\right) \, \, \mathrm{d}x.$$

\noindent
Hence from \eqref{F} and by abuse of notation clubbing the weight functions together, the above integral becomes

\begin{equation*}
    \begin{split}
        \mathcal{I}_{1}& = \frac{\sqrt{c}}{ a_{1}^{1/4}}\int_{0}^{\infty}e\left( \frac{2 \sqrt{Nxa_{2}}}{cb_{1}} \mp \frac{2 \sqrt{xNa_{1}}}{b_{1}c}\right) W(x)\, \, \mathrm{d}x \\
        & := \frac{\sqrt{c}}{ a_{1}^{1/4}} \mathcal{I}.
    \end{split}
\end{equation*}

\noindent
Now by repeated integration by parts we note that the above integral becomes negligibly small unless $a_{1} \asymp a_{2} \ll \frac{K^{2+\epsilon}}{b_{2}^{2}}$. Also by repeated integration by parts we see that the above integral becomes neglgibly small unless

\begin{equation*}
    \begin{split}
        & |\sqrt{a_{2}}-\sqrt{a_{1}}|\ll \frac{cb_{1}}{\sqrt{N}} K^\epsilon \\
        & \iff |a_{2}-a_{1}|\ll  K^\epsilon .
    \end{split}
\end{equation*}

\subsection{Further estimation} Then the $R_{2}$ sum becomes

\begin{equation*}
    \begin{split}
R_{2}& = K\sum_{b_{1}=1}^{\infty}\sum_{b_{2}=1}^{\infty}\frac{1}{b_{1}b_{2}}\sum_{c \sim K^\epsilon}\, \, \sum_{a_{2}=1}^{\infty}\frac{\lambda_{h}(a_{2})}{a_{2}^{3/4}}H\left( \frac{a_{2}b_{2}^{2}}{M}\right)\\ 
& \times \sum_{\substack{a_{1}\ll \frac{K^{2+\epsilon}}{b_{2}^{2}} \\ |a_{1}-a_{2}|\ll  K^\epsilon }} \frac{\lambda_{f}(a_{1})}{a_{1}^{1/4}}\, \mathcal{I} \,  C_{c}(a_{2}-a_{1})\\
& = K\sum_{b_{1}=1}^{\infty}\sum_{b_{2}=1}^{\infty}\frac{1}{b_{1}b_{2}}\sum_{c \sim K^\epsilon}\, \, \frac{\lambda_{h}(a_{2})}{a_{2}^{3/4}}\\ 
& \times \sum_{|q| \ll K^\epsilon}C_{c}(q)\sum_{a_{1}\ll \frac{K^{2+\epsilon}}{b_{2}^{2}} } \frac{\lambda_{f}(a_{1})\lambda_{h}(a_{1} +q)}{a_{1}^{1/4}(a_{1} +q)^{3/4}}H\left( \frac{(a_{1} +q)b_{2}^{2}}{M}\right)\, \mathcal{I} \,  
    \end{split}
\end{equation*}

\noindent
where

$$C_{c}(q)= \sum_{d | (c, q)}\mu\left( \frac{c}{d}\right)d ,$$ 

\noindent
is the Ramanujan sum. Now for $a_{1}=a_{2}$, i.e., for $q=0$ we will can proceed similarly as done in the Section \ref{diagonal}. So we will consider the off-diagonal sum, i.e., for $q=a_{2}-a_{1}\neq 0$. Now taking a smooth dyadic subdivision of the $a_{1}$-sum, we get the following

\begin{equation*}
    \begin{split}
R_{2}& \ll K\sum_{c \sim K^{\epsilon}}\sum_{0 < |q|\ll K^\epsilon}|C_{c}(q)| \left(\sup_{M^{\prime}\ll K^2 }\frac{1}{M^{\prime}}\left| \sum_{a_{1}=1 }^{\infty}\lambda_{f}(a_{1})\lambda_{h}(a_{1}+q)W\left( \frac{a_{1}}{M^{\prime}}\right
)\right| \right). 
    \end{split}
\end{equation*}

\noindent
where $W$ is some nice smooth function. Using the shifted convolution sum result (see the Theorem \ref{shifted}) we get that

\begin{equation}\label{n1}
    \begin{split}
R_{2}& \ll K\sum_{c \sim K^{\epsilon}}\sum_{0 < |q|\ll K^\epsilon}|C_{c}(q)| \left( \sup_{M^{\prime}\ll K^2 }\frac{1}{M^{\prime}} \times {M^{\prime}}^{3/4} \right)\ll K^{1/2 +\epsilon}  .
    \end{split}
\end{equation}

\noindent
This with \eqref{R1} gives the lemma \ref{E} and hence we get our Theorem \ref{main}.

\section{Shifted convolution sum estimation}
\noindent
Now we will prove the Theorem \ref{shifted}. Let us consider the shifted convolution sum

\begin{equation}\label{S}
S=\sum_{n=1}^{\infty}\lambda_{f}(n)\lambda_{g}(n+h)W\left(\frac{n}{N} \right). 
\end{equation}

 \subsection{Method of the proof of the Theorem \ref{shifted}}
We take the path of the delta method to bound the sum $S$ in \eqref{S}.  Our approach is inspired by the approach of Munshi (\cite{JEMSM}). At first, we separate the oscillations $\lambda_{f}(n)$ and $\lambda_{g}(n+h)$ using the delta symbol. Our next step is to dualize these sums.  To do that we apply the Voronoi summation formula (see \ref{Vor}). Now we use the Weil's bound for Kloosterman sum and the Ramanujan bound on average, to get the Theorem \ref{shifted}.

Here for simplification, we consider that both of the forms $f,g$ are cusp forms. Otherwise without loss of generality) if $f$ is the only cusp form then after using the Voronoi summation formula (see \ref{Vor}) to $g$ we will have a main term but for $f$, after applying the Voronoi summation formula (see \ref{Vor}) we get a savings of size $\frac{N}{Q}=\sqrt{N}$ so the size of the main term will be dominated by the bound of our Theorem \ref{shifted}.

\subsection{Preliminaries}

\subsubsection{Oscillatory integrals} Let $ \mathcal{F}$ be an be an index set and $X = X_T : \mathcal{F} \rightarrow  \mathbb{R}_{\geqslant 1}$  be a function of $ T \in \mathcal{F}. $ A family
of $ \{ \omega_T\}_{T \in \mathcal{F}}$ of smooth functions supported on a product of dyadic intervals in $ \mathbb{R}_{>0}^d$ is called $X$-inert if for each $j = (j_1, \dots, j_d ) \in \mathbb{Z}_{\geqslant 0}^d$  we have
\begin{align*}
    \sup_{T \in \mathcal{F}} \sup_{(x_1, \dots, x_d ) \in \mathbb{R}_{>0}^d } X_T^{- j_1 - \cdots - j_d} \left|  x_1^{j_1} \cdots x_d^{j_d} \omega_T^{( j_1, \cdots j_d) } ( x_1, \dots, x_d )\right| \ll_{j_1, \cdots, j_d} 1. 
\end{align*}
We study the integral of the form 

\begin{align*}
    I = \int_{\mathbb{R}} \omega (  t_1, \dots, t_d ) e (h(  t_1, \dots, t_d ))  dt_1.
\end{align*}
We have the following stationary phase lemma. 

\begin{lemma} \label{exponential integral 1}
    Suppose $ \omega= \omega_T$ is a family of $ X$-inert in $ t_1 \asymp Z$ and $ t_i \asymp X_i$ for $ i \geqslant 2$.  Suppose that on support of  $ \omega$, the function $ h $ satisfies
    \begin{align*}
        h^{ ( a_1, \dots, a_d )} ( t_1, \dots, t_d ) \ll \frac{Y}{Z^{a_1}} \frac{1}{X_2^{a_2} \cdots X_d^{a_d}}, 
    \end{align*} for all $ a_1, \cdots a_d \in \mathbb{Z}_{\geqslant 0}$.

    \begin{itemize}
        \item[$ (i)$ ] If $ h^\prime ( t_1, \dots, t_d ) \gg Y/Z$ for all $ t \in $ supp $\omega$. Suppose $Y/X \geqslant 1$. Then $ I \ll Z (Y/X)^{-A}$ for $ A$ arbitrarily large. 
        \item[$ (ii)$ ] If there exists  $ t_0 \in $ supp $\omega$ such that $ h^\prime ( t_0, \dots, t_d ) = 0$ and $ h^{\prime \prime} ( t_1, \dots, t_d ) \gg Y/Z^2$ for all $ ( t_1, \dots, t_d )\in $ supp $\omega$. Suppose $ X/Y^2 \geqslant 1$. Then

        \begin{align*}
            I = \frac{Ze^{i h( t_0, \dots, t_d )}}{ \sqrt{Y  }} W_T( t_2, \dots, t_d ) + O \left(Z  \left(\frac{Y}{X^2} \right)^{-A}\right), 
        \end{align*}
        for any $ A>0$ and for some $ X$-inert family of functions $ W_T$ supported on  $ t_0 \asymp Z$ and $ t_i \asymp X_i$ for $ i\geqslant 2$. The functions $ W_T$ depends on $A$. 
    \end{itemize}
\end{lemma}

\begin{proof}
    See \cite{BKY}.
\end{proof}

The proof of the above integral is based on the following one dimensional version, which we also require in a case where we need explicite expression for the weight function $W_T( t_2, \dots, t_d )$. 

\begin{lemma} \label{exponential integral 2}
    We consider the integral 
    \begin{align*}
        I_1 := \int_{\mathbb{R}} \omega (t) e(h(t)) dt, 
    \end{align*} where $ \omega$ is a $V$-inert function supported on the interval $ [a, b]$. Let $ 0 < \delta < 1/10$, $ Q>0$, $ Z:= Y+ b-a+1$, and assume that
    \begin{align*}
        Y \geqslant Z^{3 \delta}, \ \ \ b-a \geqslant V \geqslant \frac{Q Z^{\frac{\delta}{2}}}{ \sqrt{Y}}. 
    \end{align*}

\noindent
Suppose that there exist a unique $t_0 \in (a, b)$ such that $ h^\prime(t_0) = 0$, and furthermore

\begin{align*}
    h^{\prime \prime } (t) \gg \frac{Y}{Q^2}, \ \ \ \ h^{(j)} \ll \frac{Y}{Q^j} \ \ \  \textrm{for } \ j= 1, 2, 3, \cdots. 
\end{align*} Then we have the following asymptotic expansion for the integral
\begin{align}
    I_1 = \frac{ e(h(t_0))}{ \sqrt{h^{\prime \prime} (t_0)}} \sum_{r \leqslant 3 \delta^{-1} A} p_r(t_0) + O_{\delta, A} \left(Z^{-A} ,  \right) \ \ \ \ \ \ \textrm{where}
\end{align} 

\begin{align}
    p_r(t_0) = \frac{e^{i \pi / 4}}{ r!} \left( \frac{i}{ 4 \pi h^{\prime \prime} (t_0)}\right)^r G^{(2r)} (t_0) , \ \ \ \ \textrm{with}
\end{align}

\begin{align}
    G(t) := \omega (t) e (H(t)),  \ \ \ \ \textrm{where  } \ \ H(t) = h(t) - h(t_0) - \frac{1}{2} h^{\prime \prime} (t_0) (t - t_0)^2.
\end{align}
    
\end{lemma}

\begin{proof}
    See \cite[Section $8$]{BKY}. 
\end{proof}

\subsection{An application of the Circle-Method}
We separate the oscillations of $\lambda_{f}(n)$ and $\lambda_{g}(n)$ in the sum $S$ by using the delta symbol to arrive at
$$T = \mathop{\sum \sum}_{n, m=1}^{\infty} \lambda_{f}(n) \, \lambda_{g}(m) \,   \; W\left(\frac{n}{N}\right)\, V\left(\frac{m}{N}\right) \, \delta(n+h-m) ,$$
where $V$ is a bump function supported on the interval $[-1,1]$ and satisfies $V(0)=1$. Let $\delta : \mathbb{Z}\to \{0,1\}$ be defined by
\[
\delta(n)=
\begin{cases}
1 &\text{if}\,\,n=0 \\
0 &\text{otherwise}.
\end{cases}
\]
Then for $n\in\mathbb{Z}\cap [-2M,2M]$, we have
\begin{equation}\label{s}
\delta(n)=\frac{1}{Q}  \ \sideset{}{^\star} \sum_{a\bmod q}e\left(\frac{n a}{q}\right)\int_{\mathbb{R}}g(q,u)e\left(\frac{n u}{q	Q}\right) d u, 
\end{equation} where $Q=2M^{1/2}$. The function $g$ satisfies the following property (see $(20.158)$ and $(20.159)$ of \cite{IK}, and \cite{AG})
\begin{equation}\label{g function}
\begin{aligned}
&g(q,u)=1+h(q,u),\,\,\,\,\text{with}\,\,\,\,h(q,u)=O\left(\frac{Q}{q }\left(\frac{q}{Q}+|u|\right)^A\right),\\
&g(q,u)\ll |u|^{-A},
\end{aligned}
\end{equation}for any $A>1$.  
\begin{align} \label{deri g}
    u^j \frac{\partial^j}{ \partial u^j} g(q, u) \ll \  \min \left\lbrace \frac{Q}{q}, \frac{1}{|u|} \right\rbrace \  \log Q
\end{align}
for  $j \geqslant 1$.  In particular, the second property implies that the effective range of integral in \eqref{s} is $[-M^{\epsilon},M^{\epsilon}]$. It also follows that if $q \ll Q^{1- \epsilon}$ and $  u \ll Q^{- \epsilon} $, then  $ g(q, u)$ can
be replaced by $1$ at a cost of a negligible error term. If $ q \gg Q^{1- \epsilon}$, then we get $ u^j \frac{\partial^j}{ \partial u^j} g(q, u) \ll Q^{\epsilon}$, for any $ j \geqslant 1$. If $ q \ll  Q^{1- \epsilon}$  and $  Q^{- \epsilon} \ll |u| \ll   Q^{ \epsilon}$, then $ u^j \frac{\partial^j}{ \partial u^j} g(q, u) \ll Q^{\epsilon}$, for any $ j \geqslant 1$. Hence in all cases, we can view $  g(q, u)$ as a nice weight function. Using Parseval's identity and  Cauchy inequality, we deduce the following lemma. 

\begin{lemma} \label{L1 L2 bound  lemma}
Let $g(q, u)$ be as given in equation \eqref{s}. We have
\begin{align} \label{L1 L2 bound}
    \int_{\mathbb{R}} \left( |g (q, u)|+ |g (q, u)|^2 \right) d u \ll Q^\epsilon. 
\end{align} 
\end{lemma}
\begin{proof}
    See \cite{AG}. 
\end{proof}

 Now we let $Q =\sqrt{N}$ to write DFI's expression for $\delta(m-n)$ where $K$ will be determined later. We obtain

\begin{equation*}
    \begin{split}
      S =& \int_{\mathbb{R}}  \sum_{1 \leq q \leq Q} \frac{g(q, u)}{qQ} \sideset{}{^\star}{\sum}_{a \, \mathrm{mod} \, q} e\left(\frac{ah}{q}\right)e\left( \frac{hu}{qQ}\right)\\
      & \times \mathop{\sum }_{n=1}^{\infty} \, \lambda_{f}(n) \, e\left(\frac{an}{q}\right)e\left( \frac{nu}{qQ}\right)\,W \left(\frac{n}{N}\right) \\
      & \times \mathop{\sum }_{m=1}^{\infty} \, \lambda_{g}(m) \, e\left(-\frac{am}{q}\right)\,V\left(\frac{m}{N}\right)  e\left( -\frac{mu}{qQ}\right) \mathrm{d}u,  
    \end{split}
\end{equation*}
Trivial estimation at this stage yields $N^2$. So our job is to save $N^{1+..}$.

\subsection{Applications of Voronoi summation formula} We will apply Voronoi summation formula on the two sums:

\noindent
\textbf{\underline{For $n$:}}
Now applying the Voronoi summation formula \ref{Vor} to the $ n$-sum, we get

\begin{align*}
    S_{1} & := \mathop{\sum }_{n=1}^{\infty} \, \lambda_{f}(n) \, e\left(-\frac{an}{q}\right)\,W \left(\frac{n}{N}\right) \ e\left( \frac{nu}{qQ}\right)  \\ 
    & = \frac{1}{q}\mathop{\sum }_{n=1}^{\infty} \, \lambda_{f}(n) \, e\left(\frac{\overline{a}n}{q}\right)\, \int_0^\infty W \left(\frac{x}{N}\right)  e\left(\frac{-ux}{qQ}\right) J_{k-1}\left(\frac{4\pi\sqrt{nx}}{q}\right)\, \mathrm{d}x  \\
  & =  \frac{N}{q}\sum_{n=1}^\infty \lambda_f(n)\;e\left(\frac{\bar{a}n}{q}\right)\;\int_0^\infty  W\left(x\right) e\left(\frac{-uNx}{qQ}\right) J_{k-1}\left(\frac{4\pi\sqrt{nNx}}{q}\right) \, \mathrm{d}x.
\end{align*}

 \noindent
 Extracting the oscillation from the Bessel function we get that
\begin{align*}
     S_{1}=\frac{N^{3/4}}{\sqrt{q}}\sum_{n=1}^\infty \frac{\lambda_f(n)}{n^{1/4}}\;e\left(-\frac{\bar{a}n}{q}\right)\;\int_0^\infty e\left(\frac{uNx}{qQ}\pm \frac{2\sqrt{nNx}}{q}\right) \frac{W\left(x\right)}{x^{1/4}}\, \mathrm{d}x.
 \end{align*}

\noindent
By applying the Lemma \ref{exponential integral 1} we see that the integral is negligibly small unless 

\begin{align*}
    n \asymp \frac{Q^2}{N}= N^\epsilon .
\end{align*}

\noindent
\textbf{\underline{For $m$:}}
Applying the Voronoi summation formula \ref{Vor} to the $m$-sum, we get 

\begin{align*}
    S_{2} & := \mathop{\sum }_{m=1}^{\infty} \, \lambda_{g}(m) \, e\left(-\frac{am}{q}\right)\,V\left(\frac{m}{N}\right)  e\left( -\frac{mu}{qQ}\right) \\
    & =\frac{1}{q}\mathop{\sum }_{m=1}^{\infty} \, \lambda_{g}(m) \, e\left(\frac{\overline{a}m}{q}\right)\, \int_0^\infty e\left(\frac{uy}{qQ}\right)   \, J_{k-1}\left(\frac{4\pi\sqrt{my}}{q}\right) V\left(\frac{y}{N}\right) \, \mathrm{d}y \\ 
   &  =\frac{N^{3/4}}{\sqrt{q}}\sum_{m=1}^\infty \frac{\lambda_g(m)}{m^{1/4}}\;e\left(\frac{\bar{a}m}{q}\right)\;\int_0^\infty e\left(\frac{uNy}{qQ}\pm \frac{2\sqrt{mNy}}{q}\right) \frac{V\left(y\right)}{y^{1/4}}\, \mathrm{d}y.
\end{align*}

\noindent
By applying the Lemma \ref{exponential integral 1} we note that the above integral is negligibly small unless

\begin{align} \label{first m size}
    m \asymp \max  \frac{Q^2}{N}= N^\epsilon .
\end{align}

\noindent
We obtain the following expression of $ S$

\begin{align} \label{S after voronoi}
    S =  \frac{N^{3/2}}{Q}  \sum_{1 \leq q \leq Q} \frac{1}{q^2}  \sum_{1\leq n \ll N^\epsilon}\frac{\lambda_f(n)}{n^{1/4}}  \sum_{1\leq m\ll N^\epsilon}\frac{\lambda_g(m)}{m^{1/4}} S(h,m-n;q) \mathcal{I}(q,n,m),
\end{align} where $  S(h,m-n;q)$ is the Kloosterman sum and $ \mathcal{I}(q,n,m)$ is an integral transform which is given by

\begin{align} \label{general integral transform}
    \mathcal{I}(q,n,m) := \int_{\mathbb{R}} e\left( \frac{hu}{qQ}\right)g (u, q)  \int_0^\infty \int_0^\infty W\left(x\right) &   V\left(y\right) e\left(\frac{-u Nx}{qQ}\pm \frac{2\sqrt{nNx}}{q}\right)\, \notag\\ 
    & e\left(\frac{u Ny}{q Q}\pm \frac{2\sqrt{m Ny}}{q}\right)\, \mathrm{d}x \, \mathrm{d}y \,\mathrm{d}u.  
\end{align}

\subsection{ Simplification of the integral} We now simplify the integral transform given in equation \eqref{general integral transform}. We consider the $ u$ integral. We have

\begin{align*}
    \mathcal{I}(q,n,m) = \int_0^\infty \int_0^\infty W\left(x\right) &   V\left(y\right) e\left( \pm \frac{2\sqrt{m Ny}}{q}  \pm \frac{2\sqrt{nNx}}{q}\right)\, \notag\\ 
    &   \int_{\mathbb{R}} g (u, q)  e\left(\frac{u N (y-x+h)}{q Q}\right)\,  \mathrm{d}u \ \mathrm{d}x \ \mathrm{d}y.   
\end{align*}
Using the properties of $ g$ function from equations \eqref{g function} and \eqref{deri g}, we obtain that the integral is negligibly small unless
\begin{align*}
    |y - x| \ll |y - x+h|\ll N^\epsilon \frac{q Q}{ N} \Leftrightarrow |y - x| \ll N^\epsilon \frac{q }{Q }. 
\end{align*} Let $ x - y+h = t$. We obtain that

\begin{align*}
    \mathcal{I}(q,n,m)  & = \int_{|t-h| \ll \frac{q }{Q }}  \int_0^\infty    V\left(y\right) e\left( \pm \frac{2\sqrt{m Ny}}{q}  \pm \frac{2\sqrt{nN(y + t-h)}}{q}\right)\,
      \mathrm{d}y \,  \mathrm{d}t  \\ 
      & \hspace{2cm}+ O_A \left( N^{-A}\right),   
\end{align*}
for any positive constant $A$. Using the binomial expansion we have
\begin{align*}
    \sqrt{y + t-h} = \sqrt{y } \left( 1 + \frac{t-h}{y}\right)^{1/2}= \sqrt{y } \left( 1 + \frac{t-h}{2 y} + 
 O \left( \frac{(t-h)^2}{ y^2}\right) \right) 
\end{align*}

Using series expansion, we observe that the phase function $ e (\sqrt{nN (t-h)}/ q) $ is non-oscillating due to  size of $ t$ and $ n$, as
\begin{align*}
    \frac{\sqrt{nN}  }{q } \frac{ t-h }{ 2 \sqrt{y}} \ll N^\epsilon \frac{ Q }{q} \times  \frac{q }{Q } \ll N^\epsilon.
\end{align*}

 Hence we can absorb this into the weight function $ V$. We will consider individual bound with respect to $t$. The integral transform $ \mathcal{I}(q,n,m)$ essentially reduces to 
\begin{align*}
    \mathcal{I}(q,n,m) = \frac{q}{Q }\int_0^\infty    V\left(y\right) e\left( \pm \frac{2\sqrt{m Ny}}{q}  \pm \frac{2\sqrt{nNy }}{q}  \right)\,
      \mathrm{d}y  . 
\end{align*} There are four terms in above integral transform. We will consider one of them. Other cases are exactly same. We obtain that 
\begin{align}\label{estimate for I}
    \mathcal{I}(q,n,m) = \frac{q}{Q }\int_0^\infty    V\left(y\right) e\left( -\frac{2\sqrt{m Ny}}{q}  +  \frac{2\sqrt{nNy }}{q}  \right)\, 
      \mathrm{d}y \ll \frac{q}{Q } N^\epsilon.
\end{align}

\subsection{Further esrimates}
\noindent
Hence we have to take care of the diagonal and the off-diagonal parts.

\subsubsection{Diagonal part} Consider the case when $m=n$. For this case from \eqref{S after voronoi} we have

\begin{equation*}
    \begin{split}
         S_d &:=  \frac{N^{3/2}}{Q}  \sum_{1 \leq q \leq Q} \frac{1}{q^2}  \sum_{1\leq n \ll N^\epsilon}\frac{\lambda_f(n)\lambda_g(n)}{n^{1/2}}   C_{q}(h) \; \mathcal{I}(q,n,n),
    \end{split}
\end{equation*}

\noindent
where $C_{q}(h)$ is the Ramanujan sum. Then using \eqref{estimate for I} we get that

\begin{equation}\label{after diagonal}
    \begin{split}
         S_d &\ll  \frac{N^{3/2}}{Q}  \sum_{1 \leq q \leq Q} \frac{1}{q^2}  \times \frac{q}{Q} \ll N^{1/2 +2\theta}.
    \end{split}
\end{equation}

\subsubsection{Off-diagonal part} Consider the case when $m=n$. For this case from \eqref{S after voronoi} we have

\begin{equation*}
    \begin{split}
         S_{od} &:=  \frac{N^{3/2}}{Q}  \sum_{1 \leq q \leq Q} \frac{1}{q^2}  \sum_{1\leq n \ll N^\epsilon}\frac{\lambda_f(n)}{n^{1/4}}  \sum_{1\leq m\ll N^\epsilon}\frac{\lambda_g(m)}{m^{1/4}} S(h,m-n;q) \mathcal{I}(q,n,m) .
    \end{split}
\end{equation*}

\noindent
 Then using \eqref{estimate for I} we get that

\begin{equation}\label{after off-diagonal}
    \begin{split}
         S_{od} &\ll  \frac{N^{3/2}}{Q}  \sum_{1 \leq q \leq Q} \frac{1}{q^2}  \times \frac{q}{Q} \times \sqrt{Q}\ll N^{3/4 +2\theta}.
    \end{split}
\end{equation}

\noindent
Hence from \eqref{S after voronoi}, \eqref{after diagonal} and \eqref{after off-diagonal} we get that

\begin{equation}\label{after S}
    S \ll N^{3/4}.
\end{equation}

\noindent
This will prove the Theorem \ref{shifted}.

\
{}

\end{document}